\documentclass[11pt]{article}
\usepackage{amssymb}
\usepackage{subfigure}
\usepackage{amsbsy}
\usepackage{amsmath}
\usepackage{graphicx}
\usepackage[varg]{pxfonts}
\usepackage{epstopdf}
\usepackage{color}

\renewcommand {\thefootnote}{\fnsymbol{footnote}}

\numberwithin{equation}{section} \allowdisplaybreaks

\begin{document}
\title{\Large \bf The Minkowski problem based on the $(p,q)$-mixed quermassintegrals
\thanks{Research is supported in part by the Natural Science
Foundation of China (No.11871275; No.11371194).  }}
\author{ \small \bf Bin Chen$^{1}$\thanks{E-mail: chenb121223@163.com.}, Weidong Wang$^{2}$\thanks{E-mail: wangwd722@163.com.}, and \bf Peibiao Zhao$^{1}$\thanks{Corresponding
author E-mail: pbzhao@njust.edu.cn.}\\
\\
\small      (1. Department of Mathematics, Nanjing
University of Science and Technology, Nanjing, China) \\  \small  (2. Three Gorges Mathematical Research Center, China Three Gorges
University, Yichang, China)
}
\date{}
\maketitle
\renewcommand{\thefootnote}{\fnsymbol{footnote}}

\vskip 20pt

\begin{center}
\begin{minipage}{12.8cm}
\small
 {\bf Abstract:}  Lutwak, Yang and Zhang  \cite{L4}  introduced the concept of $L_p$ dual
 curvature measure for convex bodies and star bodies, and studied the Minkowski problem.  We  in this paper establish a new unified concept,
  in briefly, the $(p,q)$-mixed quermassintegrals, via $(p,q)$-dual mixed curvature
 measure, and further have a deep discussion on Minkowski problem with respect to the $(p,q)$-dual mixed curvature
 measure.

  By the way, we derive at some important properties and geometric inequalities for $(p,q)$-mixed quermassintegrals.



 \vskip 0.1cm
\end{minipage}
\end{center}

\vskip 20pt
\section{\bf Introduction}
~~~~Brunn and Minkowski pioneered the classical Brunn-Minkowski theory of convex bodies in $n$-dimensional Euclidean space. The Minkowski linear combination of vectors, mixed volumes, surface area measures and basic Brunn-Minkowski inequalities etc. are the essence in Brunn-Minkowski theory. The classical dual Brunn-Minkowski theory of star bodies was introduced by Lutwak (\cite{L1}) in 1975. The formation of dual Brunn-Minkowski theory was based on radial Minkowski combination and dual mixed volumes of star bodies. For more information about classical Brunn-Minkowski theory and dual Brunn-Minkowski theory, please refer to two excellent books \cite{G1,S0} in details.

In modern convex geometry, the $L_p$ Brunn-Minkowski theory and $L_p$ dual Brunn-Minkowski theory generalize and dualize the classical Brunn-Minkowski theory. The $L_p$ mixed quermassintegrals $W_{p,i}(M,N)$ of convex bodies introduced by Lutwak (see\cite{L2}) and $L_p$ dual mixed quermassintegrals $\widetilde{W}_{-p,i}(M,N)$ of star bodies introduced by Wang and Leng (see\cite{W}) generalized the $L_p$ mixed volumes $V_p(M,N)$ introduced in \cite{L2} and $L_p$ dual mixed volumes $\widetilde{V}_{-p}(M,N)$ introduced in \cite{L3} are the main geometric measures in these two theories. In 2018, Lutwak, Young and Zhang established the $(p,q)$-mixed volume $V_{p,q}(M,N,Q)$ in \cite{L4}. Meanwhile, they demonstrated a surprising connection between the $L_p$ mixed volume and $L_p$ dual mixed volume.
During the last three decades, scholars have done a lot of work on the research of $L_p$ Brunn-Minkowski theory and $L_p$ dual Brunn-Minkowski theory, see e.g., \cite{BF,BL,BL1,C0,FW,H,HL1,L0,LWZ,W0,WH,WL,WW,ZZ}.

The Minkowski problem in Brunn-Minkowski theory and dual Brunn-Minkowski theory is a hot topic of research.  The classical Minkowski problem was first studied by Minkowski \cite{M1,M2}, who proved both existence and uniqueness of solutions for the given measure. The $L_p$ Minkowski problem in $L_p$ Brunn-Minkowski theory was introduced in \cite{L2}. Existence of solutions for the dual Minkowski problem for even data within the class of origin-symmetric convex bodies was proved in \cite{H}. Very recently, $L_p$ dual Minkowski problem for $L_p$ dual curvature measures was stated in \cite{L4}, which a unified Minkowski problem of $L_p$ Minkowski problem when $q=n$ and dual Minkowski problem when $p=0$ and $N=B$. The Minkowski problem and dual Minkowski problem have already aroused considerable attention; see, e.g., \cite{Ch,C0,CL,HS,HL,LO,LY,J,WF,Z3,Z4,Z5,ZY,ZY1}.

The $L_p$ dual curvature measure $\widetilde{C}_{p,q}(M,N,\cdot)$ of $M\in\mathcal{K}_o^n$ and $N\in\mathcal{S}_o^n$ was introduced in \cite{L4}, which constructs a unified frame for $L_p$ surface area measures introduced in \cite{L2} and dual curvature measures introduced in \cite{H}. Meanwhile, they demonstrated a surprising connection between the $L_p$ Brunn-Minkowski theory and $L_p$ dual Brunn-Minkowski theory by establishing geometric inequalities and variational integral formulas of $L_p$ dual curvature measure, respectively.

Based on the concepts introduced above, it is worth researching to find a concept
that can unify the $L_p$ mixed quermassintegrals, $L_p$ dual mixed quermassintegrals and $(p,q)$-mixed volume.
Motivated by ideas in \cite{L2} and \cite{L4}, the main purpose of
this paper is to introduce the concept of $(p,q)$-dual mixed
curvature measures of convex bodies and star bodies. Next, by
the concept of $(p,q)$-dual mixed curvature measures,
 we establish the integral formula of $(p,q)$-mixed quermassintegrals
$\widetilde{W}_{p,q,i}(M,N,Q)$, which is a unified structure  of
$(p,q)$-mixed volumes, $L_p$ mixed quermassintegrals and $L_p$ dual
mixed quermassintegrals.
Our main work in this article can be illustrated in the block diagram below:

\begin{figure}[h!]
\centering
\begin{minipage}{15cm}
\centering
\includegraphics[scale=0.7]{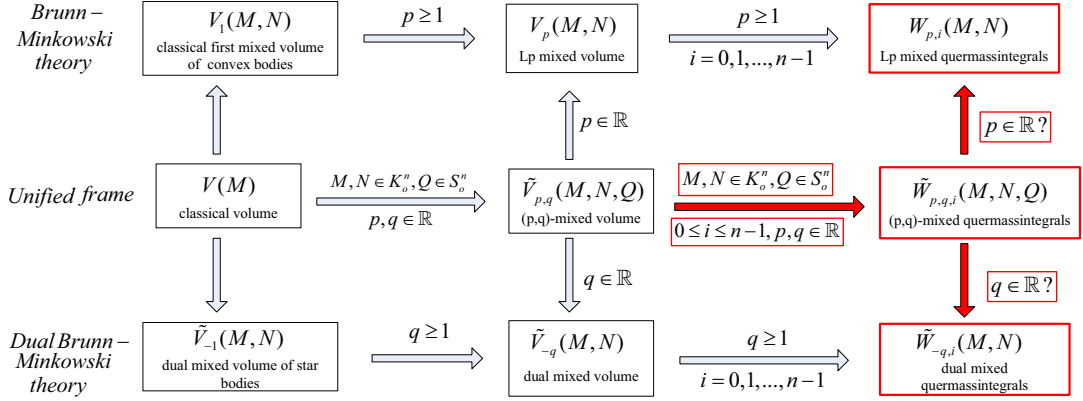}
\end{minipage}
\caption{Idea guide chart.}
\end{figure}

Let $M$ be a convex body if $M$ is a compact, convex subset in $n$-dimensional Euclidean space $\mathbb{R}^n$ with non-empty interior. The set of all convex bodies containing the origin in their interiors in $\mathbb{R}^n$ is written as $\mathcal{K}_o^n$. We write $u$ for a unit vector and $B$ for the unit ball centered at the origin, the surface of $B$ denoted by $S^{n-1}$. We shall use $V(M)$ for the $n$-dimensional volume of the body $M$ in $\mathbb{R}^n$, we write $V(B)=\omega_n$.

For $M\in\mathcal{K}^n$, the support function $h_M=h(M,\cdot): \mathbb{R}^n\rightarrow\mathbb{R}$ is defined by
$$h(M,x)=\max\{x\cdot y: y\in M\},\ \ \ \ x\in\mathbb{R}^n.$$
where $x\cdot y$ denotes the standard inner product of $x$ and $y$ in $\mathbb{R}^{n}$ (see\cite{G1,S0}).

In the dual Brunn-Minkowski theory, for  $M$ is a compact star shaped (about the origin) in $\mathbb{R}^n$, the radial function, $\rho_M=\rho(M,\cdot)$: $\mathbb{R}^n\backslash\{0\}\rightarrow[0,+\infty)$, is defined by
  $$\rho(M,x)=\max\{\lambda\geq0:\lambda x\in M\},\ \ \ \   x\in \mathbb{R}^n\backslash\{0\}.$$
If $\rho_M$ is positive and continuous, then call $M$ is a star body and write as $\mathcal{S}_o^n$ (see\cite {G1,S0}).

For $K,L\in\mathcal{K}_o^n$, $p\geq1$ and $\lambda,\mu\geq0$ (not
both zero), the $L_p$-Minkowski linear  combination, $\lambda\cdot
K+_p\mu\cdot L\in\mathcal{K}_o^n$, of $K$ and $L$ satisfies
(see\cite{F})
$$h(\lambda\cdot K+_p\mu\cdot L,\cdot)^p=\lambda h(K,\cdot)^p+\mu h(L,\cdot)^p.$$

The $L_p$ mixed surface area measures $S_{p,i}(M,\cdot)$ of $M\in\mathcal{K}_o^n$ can be defined by the variational formula,
$$\frac{d}{dt}W_{p,i}(M+_pt\cdot N)\bigg|_{t=0^+}=\frac{1}{p}\int_{S^{n-1}}h_N^p(u)dS_{p,i}(M,u),$$
for $N\in\mathcal{K}_o^n$ and $0\leq i\leq n-1$.

For $q\in \mathbb{R}$ and $0\leq i\leq n-1$, the $q$-th dual mixed quermassintegral of $M,N\in\mathcal{S}_o^n$, is defined by
\begin{align}\label{1.1}
\widetilde{W}_{q,i}(M,N)=\frac{1}{n}\int_{S^{n-1}}\rho_M^q(u)\rho_N^{n-q-i}(u)du.
\end{align}

In this paper, the $(p,q)$-dual mixed curvature measures, $\widetilde{\mathcal{C}}_{p,q,j}$, is a three-parameter family of Borel measures on $S^{n-1}$. For $p,q\in\mathbb{R}$, $j\neq n$, and $M\in\mathcal{K}_o^n$, $Q\in\mathcal{S}_o^n$, we define the $\widetilde{\mathcal{C}}_{p,q,j}(M,Q,\cdot)$ by
$$\int_{S^{n-1}}g(v)d\widetilde{\mathcal{C}}_{p,q,j}(M,Q,v)
=\frac{1}{n}\int_{S^{n-1}}g(\alpha_M(u))h_M^{-p}(\alpha_M(u))\rho_M^q(u)\rho_Q^{n-q-j}(u)du,$$
for a continuous function $g: S^{n-1}\rightarrow\mathbb{R}$, where $\alpha_M$ is the radial Gauss map (see Section 2) that associates to almost $u\in S^{n-1}$ the unique out unit normal vector at the point $\rho_M(u)u\in\partial M$.

The $L_p$ mixed surface area measures and $L_p$ dual curvature measures are special cases of the $(p,q)$-dual mixed curvature measures. For $p,q\in\mathbb{R}$, $M\in\mathcal{K}_o^n$ and $M\in\mathcal{S}_o^n$
$$\widetilde{\mathcal{C}}_{p,q,0}(M,Q,\cdot)=\widetilde{\mathcal{C}}_{p,q}(M,Q,\cdot),\ \ \ \ \widetilde{\mathcal{C}}_{p,q,0}(M,M,\cdot)=\frac{1}{n}S_{p}(M,\cdot),$$
$$\widetilde{\mathcal{C}}_{p,n,0}(M,B,\cdot)=\frac{1}{n}S_{p}(M,\cdot),\ \ \ \ \widetilde{\mathcal{C}}_{0,q,0}(M,B,\cdot)=\frac{1}{n}\widetilde{\mathcal{C}}_{q}(M,\cdot).$$

Recall that $L_p$ Minkowski linear combination, we give the variational formulas of $q$-th dual mixed quermassintegral: Suppose $p,q\in\mathbb{R}$, $j\neq n$. For $M,N\in\mathcal{K}_o^n$ and $Q\in\mathcal{S}_o^n$,
\begin{align}\label{1.2}
\frac{d}{dt}\widetilde{W}_{q,j}(M+_pt\cdot N,Q)\bigg|_{t=0^+}
=\frac{q}{p}\int_{S^{n-1}}h_N^{p}(u)d\widetilde{C}_{p,q,j}(M,Q,u).
\end{align}

It will be shown that the $(p,q)$-mixed quermassintegral $\widetilde{W}_{p,q,j}$ has a integral representation using the $(p,q)$-dual mixed curvature measure: For $p,q\in\mathbb{R}$, and $j=0,1,...,n-1$, there exists a regular Borel measure $\mathcal{\widetilde{C}}_{p,q,j}(M,N,\cdot)$ on $S^{n-1}$, such that the $(p,q)$-mixed quermassintegral $\widetilde{W}_{p,q,j}$ has the following integral representation:
\begin{align}\label{1.3}
\widetilde{W}_{p,q,j}(M,N,Q)
=\int_{S^{n-1}}h_N^p(v)d\widetilde{C}_{p,q,j}(M,Q,v).
\end{align}

The $L_p$ mixed quermassintegrals and the $L_p$ dual mixed quermassintegrals will be the special cases of $(p,q)$-mixed quermassintegrals,
\begin{align}\label{1.4}
\widetilde{W}_{p,q,j}(M,N,M)=W_{p,j}(M,N),
\end{align}
\begin{align}\label{1.5}
\widetilde{W}_{p,q,j}(M,M,Q)=\widetilde{W}_{q,j}(M,N).
\end{align}

The following inequality for $(p,q)$-mixed quermassintegral is a generalization of the $L_p$ Minkowski  inequality: Suppose $p,q\in\mathbb{R}$ satisfies $1\leq\frac{q}{n-j}\leq p$, $j\neq n$. If $M,N\in\mathcal{K}_o^n$ and $Q\in\mathcal{S}_o^n$, then
$$\widetilde{W}_{p,q,j}(M,N,Q)\geq W_j(M)^{\frac{q-p}{n-j}}W_j(N)^{\frac{p}{n-j}}\widetilde{W}_j(Q)^{\frac{n-q-j}{n-j}},$$
with equality if and only if $M,N,Q$ are dilates when $1<\frac{q}{n-j}<p$, while only $M$ and $N$ are dilates when $q=n$ and $p>1$, and $M$ and $N$ are homotheic when $q=n$ and $p=1$.

Therefore, we can see that the $L_p$ Brunn-Minkowski theory and $L_p$ dual Brunn-Minkowski theory can be finally unified by definitions of $(p,q)$-dual mixed curvature measures and $(p,q)$-mixed quermassintegrals.

The another goal of this article is to study the Minkowski problem for $(p,q)$-dual mixed curvature measures, which is a general Minkowski problem that unify the $L_p$ Minkowski problem and the $L_p$ dual Minkowski problem. To prove the existence of Minkowski problems, the key is to convert the Minkowski problem to
maximization problem.

The existence problem for $(p,q)$-dual mixed curvature measures is: {\it For $p,q\in\mathbb{R}$, $j\neq n$ and $Q\in\mathcal{S}_o^n$. Given a Borel measure $\mu$, what are necessary and sufficient conditions on $\mu$ so that there exists a convex body $M\in\mathcal{K}_o^n$ whose dual curvature measure $\widetilde{\mathcal{C}}_{p,q,j}(M,Q,\cdot)$ is the given measure $\mu$?}

The uniqueness problem for $(p,q)$-dual mixed curvature measures is: {\it For $p,q\in\mathbb{R}$, $j\neq n$ and $Q\in\mathcal{S}_o^n$. If $M_1,M_2\in\mathcal{K}_o^n$, are such that
$$\widetilde{\mathcal{C}}_{p,q,j}(M_1,Q,\cdot)=\widetilde{\mathcal{C}}_{p,q,j}(M_2,Q,\cdot),$$
then how is $M_1$ related to $M_2$?}

The organization is as follows. The corresponding background
materials and some results are introduced in Section 2. In Section 3, we define the concept of $(p,q)$-dual mixed curvature measures on the basis of $L_p$ dual curvature measures. In Section 4, we then establish the variational formulas for $q$-th dual mixed quermassintegrals, and obtain the integral formula for $(p,q)$-mixed quermassintegrals via the concept of $(p,q)$-dual mixed curvature measures. In Section 5, we derive at some important geometric inequalities involving Minkowski inequality etc. Finally, we further study the Minkowski problem for $(p,q)$-dual mixed curvature measures in Section 6.

\section{\bf Preliminaries}
$\quad$In this section, we give the interrelated background materials and some results. Gardner's book \cite{G1} and Schneider's book \cite{S0} are our standard references for the basics regarding convex bodies and star bodies.

{\bf 2.1~~Support function, radial function, Wulff shape and convex hull}

For $M\in\mathcal{K}^n$, the support function $h_M=h(M,\cdot): \mathbb{R}^n\rightarrow\mathbb{R}$ is given by
$$h(M,x)=\max\{x\cdot y: y\in M\},\ \ \ \ x\in\mathbb{R}^n.$$
where $x\cdot y$ denotes the standard inner product of $x$ and $y$ in $\mathbb{R}^{n}$. For $\phi\in GL(n)$ (which denotes the general linear transformation group), the image of $M$ under $\phi$, we have
$$h(\phi M,x)=h(M,\phi^tx),\ \ \ \ x\in\mathbb{R}^n,$$
where $\phi^t$ denotes the transpose of $\phi$.

For $M\in\mathcal{S}_o^n$, the radial function, $\rho_M=\rho(M,\cdot)$: $\mathbb{R}^n\backslash\{0\}\rightarrow[0,+\infty)$, is defined by
  $$\rho(M,x)=\max\{\lambda\geq0:\lambda x\in M\},\ \ \ \   x\in \mathbb{R}^n\backslash\{0\}.$$
Similarly, for $\phi\in GL(n)$, we have
$$\rho(\phi M,x)=\rho(M,\phi^{-1}x),\ \ \ \ x\in\mathbb{R}^n\backslash\{0\}.$$

For $M\in\mathcal{K}_o^n$, the polar body $M^\ast$ of $M$ is the convex body in $\mathbb{R}^n$, and defined by
$$M^\ast=\{x\in\mathbb{R}^n: x\cdot y\leq1,\ for\ all\ y\in M\}.$$

From this definition, we easily see that on $\mathbb{R}^n\backslash\{0\}$,
\begin{align}\label{2.1}
\rho_M=\frac{1}{h_{M^\ast}}\ \ \ \ and \ \ \ \ h_M=\frac{1}{\rho_{M^\ast}}.
\end{align}
It follows that for $M\in\mathcal{K}_o^n$
$$(M^\ast)^\ast=M.$$

Let $\Omega\subset S^{n-1}$ denote a set that is closed and cannot be contained in any closed hemisphere of $S^{n-1}$. The Wulff shape $[h]\in\mathcal{K}_o^n$, of a continuous function $h: \Omega\rightarrow(0,\infty)$, also known as the Aleksandrov body of $h$, is define by
$$[h]=\bigcap_{u\in\Omega}\{x\in\mathbb{R}^n: x\cdot u\leq h(u)\}.$$
Because of the restrictions placed on $\Omega$, we see that $[h]\in\mathcal{K}_o^n$. If $M\in\mathcal{K}_o^n$, then
$$[h_M]=M.$$

For the radial function of the Wulff shape, we have
$$\rho_{[h]}(u)^{-1}=\max_{u\in\Omega}(u\cdot v)h(v)^{-1}.$$

Let $\rho: \Omega\rightarrow(0,\infty)$ be continuous and $\{\rho(u)u: u\in S^{n-1}\}$ is a compact set in $\mathbb{R}^n$. The convex hull $\langle\rho\rangle$ is defined by
$$\langle\rho\rangle=conv\{\rho(u)u: u\in S^{n-1}\},$$
is compact as well (see Schneider \cite{S0} Theorem 1.1.11). If $M\in\mathcal{K}_o^n$, then
$$\langle\rho_M\rangle=M.$$

For the support function of the convex hulls, we have
$$h_{\langle\rho\rangle}(v)\max_{u\in\Omega}(u\cdot v)\rho(u),$$
for all $v\in S^{n-1}$.

Using the concept of Wulff shape, the definition of $L_p$ Minkowski combination can be extended to $p<1$ and even negative $k$ or $l$: Fix a real $p\neq0$. For $M,N\in\mathcal{K}_o^n$, and $k,l\in\mathbb{R}$ such that $kh_M^p+lh_N^p$ is a strictly positive function on $S^{n-1}$, define the $L_p$ Minkowski combination, $k\cdot M+_pl\cdot N\in\mathcal{K}_o^n$, by
\begin{align}\label{2.2}
k\cdot M+_pl\cdot N=[(kh_M^p+lh_N^p)^{\frac{1}{p}}].
\end{align}
For $\phi\in GL(n)$ and $p\neq0$,
\begin{align}\label{2.3}
k\cdot\phi M+_pl\cdot\phi N=\phi(k\cdot M+_pl\cdot N).
\end{align}

{\bf 2.2~~The radial Gauss map}

Suppose $M\in\mathcal{K}^n$ in $\mathbb{R}^n$. For $v\in\mathbb{R}^n\backslash\{0\}$, the hyperplane
$$H_M(v)=\{x\in\mathbb{R}^n: x\cdot v=h_M(v)\}$$
is called the supporting hyperplane to $M$ with outer normal $v$.

The spherical image of $\sigma\subset\partial M$ is defined by
$$\mathbb{V}_M(\sigma)=\{v\in S^{n-1}: x\in H_M(v)\ for \ some \ x\in\sigma\}\subset S^{n-1}.$$
The reverse spherical image of $\eta\subset S^{n-1}$ is defined by
$$\mathbb{X}_M(\eta)=\{x\in\partial M: x\in H_M(v) \ for \ some \ v\in\eta\}\subset\partial M.$$

Let $\sigma_M\subset\partial M$ be the set consisting of all $x\in \partial M$, for which the set $\mathbb{V}_M(\{x\})$, which we frequently abbreviate as $\mathbb{V}_M(x)$, contains more than a single element. It is well known that $\mathcal{H}^{n-1}(\sigma_M)=0$ (see p. 84 of Schneider \cite{S0}). On precisely the set of regular radial vectors of $\partial M$ is defined the function
$$\nu_M:\partial M\backslash\sigma_M\rightarrow S^{n-1},$$
by letting $\nu_M(x)$ be the unique element in $\mathbb{V}_M(x)$ for each $x\in\partial M\backslash\sigma_M$. The function $\nu_M$ is called the spherical image map of $M$ and is known to be continuous (see Lemma 2.2.12 of Schneider \cite{S0}). It will occasionally be convenient to abbreviate $\partial M\backslash\sigma_M$ by $\partial^{\prime}M$. Since $\mathcal{H}^{n-1}(\sigma_M)=0$, when the integration is with respect to $\mathcal{H}^{n-1}$, it will be immaterial if the domain is over subsets of $\partial^{\prime} M$ or $\partial M$.

For $M\in\mathcal{K}_o^n$ and $u\in S^{n-1}$, the radial map of $M$ is defined as,
$$r_M: S^{n-1}\rightarrow\partial M\ \ by\ \ r_M(u)=\rho_M(u)u\in\partial M.$$

For $\omega\subset S^{n-1}$, the radial Gauss image of $\omega$ is defined as
$$\Re_M(\omega)=\{v\in S^{n-1}: r_M(u)\in H_M(v)\ for\ some\ u\in\omega\},$$
and thus, for $u\in S^{n-1}$,
\begin{align}\label{2.4}
\Re_M(u)=\{v\in S^{n-1}: r_M(u)\in H_M(v)\}.
\end{align}

For $M\in\mathcal{K}_o^n$, the radial Gauss map is defined by
$$\alpha_M: S^{n-1}\backslash\omega_M\rightarrow S^{n-1},\ \ by \ \ \alpha_M=\nu_M\circ r_M$$
where $\omega_M=r_M^{-1}(\sigma_M)$. Since $r_M^{-1}$ is a bi-Lipschitz map between the space $\partial M$ and $S^{n-1}$ it follows that $\omega_M$ has spherical Lebesgue measure zero. Observe that if $u\in S^{n-1}\backslash\omega_M$, then $\Re_M(u)$ contains only the element $\alpha_M(u)$. For $x\in\mathbb{R}^n\backslash\{0\}$, define $\overline{x}\in S^{n-1}$ by $\overline{x}=x/|x|$. Since both $\nu_M$ and $r_M$ are continuous, then $\alpha_M$ is continuous. Note that for $x\in\partial^{\prime}M$,
\begin{align}\label{2.5}
\alpha_M(\overline{x})=\nu_M(x),
\end{align}
and hence, for $x\in\partial^{\prime}M$,
\begin{align}\label{2.6}
h_M(\alpha_M(\overline{x}))=h_M(\nu_M(x))=x\cdot\nu_M(x).
\end{align}

For $\eta\subset S^{n-1}$, the reverse radial Gauss image of $\eta$ is defined by
$$\Re_M^\ast(\eta)=r_M^{-1}(\mathbb{X}_M(\eta))=\langle\mathbb{X}_M(\eta)\rangle.$$
Thus,
$$\Re_M^\ast(\eta)=\{\overline{x}: x\in\partial M\ \ where\ x\in H_M(v)\ for\ some\ v\in\eta\}.$$
If $\eta$ contains only the single vector $v\in S^{n-1}$, then
$$\Re_M^\ast(v)=\{\overline{x}: x\in\partial M\ \ where\ x\in H_M(v)\}.$$

Note that the set $\eta\subset S^{n-1}$,
$$\Re_M^\ast(\eta)=\{u\in S^{n-1}: r_M(u)\in H_M(v)\ for\ some\ v\in\eta\},$$
and for $u\in S^{n-1}$ and $\eta\subset S^{n-1}$, we see that from (\ref{2.4}) that
\begin{align}\label{2.7}
u\in\Re_M^\ast(\eta)\ \ \Leftrightarrow\ \ \Re_M(u)\cap \eta\neq\varnothing.
\end{align}
If $u\not\in\omega_M$, then $\Re_M(u)=\{\alpha_M(u)\}$, then (\ref{2.7}) becomes
\begin{align}\label{2.8}
u\in\Re_M^\ast(\eta)\ \Leftrightarrow\ \alpha_M(u)\in\eta,
\end{align}
and hence (\ref{2.8}) holds for almost all $u\in S^{n-1}$, with respect to spherical Lebesgue measure.

The reverse radial Gauss image of a convex body and the radial Gauss image of its polar body are related (see\cite{H}).

\noindent{\bf Lemma 2.1}~~{\it For each $\eta\subset S^{n-1}$. If $M\in\mathcal{K}_o^n$, then
$$\Re_M^\ast(\eta)=\Re_{M^\ast}(\eta).$$

For almost all $v\in S^{n-1}$, then $\Re_M^\ast(v)=\{\alpha_M^\ast(v)\}$,
and $\Re_{M^\ast}(v)=\{\alpha_{M^\ast}(v)\}$. Therefore,
$$\alpha_M^\ast=\alpha_{M^\ast}.$$}

{\bf 2.3~~The surface area measure}

By using the spherical image and reverse spherical image, one can define the surface area measures, and their $L_p$ extensions.

For Borel set $\eta\subset S^{n-1}$ and $M\in\mathcal{K}_o^n$. the surface area measures $S(M,\cdot)$ can be defined by
$$S(M,\eta)=\mathcal{H}^{n-1}(\mathbb{X}_M(\eta)).$$

For convex bodies $M,N$ in $\mathbb{R}^n$, the mixed quermassintegral $W_i(M,N)$ has the following integral representation:
$$W_i(M,N)=\frac{1}{n}\int_{S^{n-1}}h_N(u)dS_i(M,u),$$
where $dS_i(M,u)/dS(u)=f_i(M,u): S^{n-1}\rightarrow [0,\infty)$ denote the $i$-th curvature function (see\cite{L2}).

For $p\in\mathbb{R}$ and $M,N\in\mathcal{K}_o^n$, the $L_p$ mixed quermassintegral $W_{p,i}(M,N)$ is defined by
\begin{align}\label{2.10}
W_{p,i}(M,N)=\frac{1}{n}\int_{S^{n-1}}h_N^p(u)dS_{p,i}(M,u),
\end{align}
where the $S_{p,i}(M,\cdot)$ is the $L_p$ surface area measure.

For $p\geq1$, the $L_p$ Minkowski inequality of $L_p$ mixed quermassintegral is obtained
\begin{align}\label{2.11}
W_{p,i}(M,N)\geq W_i(M)^{\frac{n-p-i}{n-i}}W_i(N)^{\frac{p}{n-i}},
\end{align}
with equality if and only if $M$ and $N$ are dilates when $p>1$, and $M$ and $N$ are homothets when $p=1$.

The following important integral identity was established in \cite{H}: For $q\in\mathbb{R}$, $j\neq n$ and $M\in\mathcal{K}_o^n$. If $f: S^{n-1}\rightarrow\mathbb{R}$ is bounded and Lebesgue integrable, then
\begin{align}\label{2.12}
\int_{S^{n-1}}f(u)\rho_M^{q}(u)du
=\int_{\partial M}f(\overline{x})|x|^{q-n}(x\cdot\nu_M(x))d\mathcal{H}^{n-1}(x).
\end{align}

\section{\bf The $(p,q)$-dual mixed curvature measures}

$\quad$In this section, we introduce the concept of $(p,q)$-dual mixed curvature measures on the basis of $L_p$ dual curvature measures.

Let $Q\in\mathcal{S}_o^n$, define a continuous and positively function $\|\cdot\|_Q: \mathbb{R}^n\rightarrow[0,\infty)$ by
\begin{align}\label{3.1}
\|x\|_Q=\bigg\{_{0,\ \ \ \ \ \ \ \ \ \ \ \ \ \ \ x=0.}^{\rho_Q^{-1}(x),\ \ \ \ \ \ \ \ \ \ x\neq0}
\end{align}
Note that $\|\cdot\|_Q$ is a continuous, and a positively homogeneous function of degree one. When $Q$ is an origin-symmetric convex bodies in $\mathbb{R}^n$, then $\|\cdot\|_Q$ is an ordinary norm in $\mathbb{R}^n$, and $(\mathbb{R}^n,\|\cdot\|_Q)$ is the $n$-dimensional Banach space whose unit ball is $Q$.

\noindent{\bf Definition 3.1}~~{\it Suppose $q\in\mathbb{R}$ and $j\neq n$. For $M\in\mathcal{K}_o^n$, $Q\in\mathcal{S}_o^n$ and each Borel $\eta\subseteq S^{n-1}$, the $q$-th dual mixed curvature measure, $\widetilde{\mathcal{C}}_{q,j}(M,Q,\cdot)$, of $M,Q$ is defined by
\begin{align}\label{3.2}
\nonumber\widetilde{\mathcal{C}}_{q,j}(M,Q,\eta)&=\frac{1}{n}\int_{\Re_M^\ast(\eta)}\rho_M^q(u)\rho_Q^{n-q-j}(u)du\\
&=\frac{1}{n}\int_{S^{n-1}}\mathbb{I}_{\Re^*_M(\eta)}(u)\rho_M^q(u)\rho_Q^{n-q-j}(u)du.
\end{align}
Moreover, for any $p\in\mathbb{R}$, the $(p,q)$-th dual mixed curvature measure,  $\widetilde{\mathcal{C}}_{p,q,j}(M,Q,\cdot)$, is defined by
\begin{align}\label{3.3}
\frac{d\widetilde{\mathcal{C}}_{p,q,j}(M,Q,\cdot)}{d\widetilde{\mathcal{C}}_{q,j}(M,Q,\cdot)}=h^{-p}(M,\cdot).
\end{align}}
In particular, if $p=0$ in (\ref{3.3}), then
$$\widetilde{\mathcal{C}}_{0,q,j}(M,Q,\cdot)=\widetilde{\mathcal{C}}_{q,j}(M,Q,\cdot).$$

\noindent{\bf Remark}~~In particular, let $j=0$ in (\ref{3.2}), the dual curvature measures $\widetilde{\mathcal{C}}_{q,0}(M,N,\eta)=\widetilde{\mathcal{C}}_{q}(M,N,\eta)$. If $j=0$ in (\ref{3.3}), then $d\widetilde{\mathcal{C}}_{p,q,0}(M,N,\cdot)=h_K^{-p}d\widetilde{\mathcal{C}}_{p.q}(M,N,\cdot)$. Further, if $p=0$, then $\widetilde{\mathcal{C}}_{0,q}(M,N,\cdot)=\widetilde{\mathcal{C}}_{q}(M,N,\cdot)$ (see\cite{L4}).

Then, we are ready to establish the integral representation of $(p,q)$-dual mixed curvature measures. But we also need a tool as follows. Let $\varphi: S^{n-1}\rightarrow\mathbb{R}$ be a simple function on $S^{n-1}$ by letting
$$\varphi=\sum_{i=1}^kc_i\mathbb{I}_{\eta_i}$$
for $c_i\in\mathbb{R}$ and Borel set $\eta_i\subset S^{n-1}$.

\noindent{\bf Lemma 3.1}~~{ \it Suppose $q\in\mathbb{R}$ and $j\neq n$. For $M\in\mathcal{K}_o^n$, $Q\in\mathcal{S}_o^n$ and a bounded Borel function $f: S^{n-1}\rightarrow\mathbb{R}$, we have
\begin{align}\label{3.4}
\int_{S^{n-1}}f(v)d\widetilde{\mathcal{C}}_{q,j}(M,Q,v)
=\frac{1}{n}\int_{S^{n-1}}f(\alpha_M(v))\rho_M^q(v)\rho_Q^{n-q-j}(v)dv.
\end{align}}

{\it \bf Proof.}~~From (\ref{3.2}) and the fact (\ref{2.8}), for any Borel sets $\eta\subset S^{n-1}$, we get
\begin{align}\label{3.5}
\nonumber\int_{S^{n-1}}\varphi(u)d\widetilde{\mathcal{C}}_{q,j}(M,Q,u)&=\int_{S^{n-1}}\sum_{i=1}^kc_i\mathbb{I}_{\eta_i}(u)
d\widetilde{\mathcal{C}}_{q,j}(M,Q,u)\\
\nonumber&=\sum_{i=1}^kc_i\int_{S^{n-1}}\mathbb{I}_{\eta_i}(u)
d\widetilde{\mathcal{C}}_{q,j}(M,Q,u)\\
\nonumber&=\sum_{i=1}^kc_i\widetilde{\mathcal{C}}_{q,j}(M,Q,\eta_i)\\
&=\frac{1}{n}\sum_{i=1}^kc_i\int_{\Re_M^\ast(\eta_i)}\rho_M^q(v)\rho_Q^{n-q-j}(v)dv\\
\nonumber&=\frac{1}{n}\int_{S^{n-1}}\sum_{i=1}^kc_i\mathbb{I}_{\Re_M^\ast(\eta_i)}(v)\rho_M^q(v)\rho_Q^{n-q-j}(v)dv\\
\nonumber&=\frac{1}{n}\int_{S^{n-1}}\sum_{i=1}^kc_i\mathbb{I}_{{\eta_i}}(\alpha_M(v))\rho_M^q(v)\rho_Q^{n-q-j}(v)dv\\
\nonumber&=\frac{1}{n}\int_{S^{n-1}}\varphi(\alpha_M(u))\rho_M^q(u)\rho_Q^{n-q-j}(u)du.
\end{align}

In (\ref{3.5}), we choose a sequence of simple functions $f_k$ that converge uniformly to $f$. Then $f_k\circ\alpha_M$ converges to $f\circ\alpha_M$ a.e. with respect to spherical Lebesgue measure. Since $f$ is a Borel function on $S^{n-1}$ and the radial Gauss map $\alpha_M$ is continuous on $S^{n-1} \backslash \omega_M$, then composite function $f\circ\alpha_M$ is a Borel function on $S^{n-1}\backslash \omega_M$ as well. Since $f$ is bounded and $\omega_M$ has Lebesgue measure zero. Thus, $f$ and $S^{n-1}\backslash \omega_M$ are Lebesgue integrable on $S^{n-1}$. Taking the $k\rightarrow\infty$, we can establish (\ref{3.4}).
\hfill${\square}$

From (\ref{3.3}), (\ref{3.4}) and (\ref{2.8}), we can easily get the integral formula for $(p,q)$-dual mixed curvature measures as follows:
Suppose $p,q\in\mathbb{R}$ and $j\neq n$. For $M\in\mathcal{K}_o^n$, $Q\in\mathcal{S}_o^n$ and each Borel set $\eta\subseteq S^{n-1}$, then
\begin{align}\label{3.6}
\widetilde{\mathcal{C}}_{p,q,j}(M,Q,\eta)
=\frac{1}{n}\int_{\Re_M^\ast(\eta)}h_M^{-p}(\alpha_M(u))\rho_M^q(u)\rho_Q^{n-q-j}(u)du.
\end{align}

Next, we are going to get some properties of $(p,q)$-dual mixed curvature measures.

\noindent{\bf Lemma 3.2}~~{\it Suppose $p,q\in\mathbb{R}$ and $j\neq n$. If $M\in\mathcal{K}_o^n$ and $Q\in\mathcal{S}_o^n$, then for each bounded, Borel function $g: S^{n-1}\rightarrow\mathbb{R}$
\begin{align}\label{3.7}
\int_{S^{n-1}}g(v)d\widetilde{\mathcal{C}}_{p,q,j}(M,Q,v)
=\frac{1}{n}\int_{S^{n-1}}g(\alpha_M(u))h_M^{-p}(\alpha_M(u))\rho_M^q(u)\rho_Q^{n-q-j}(u)du.
\end{align}}
For the right hand of formula (\ref{3.7}), the integration is with respect to spherical Lebesgue measure.

{\it \bf Proof.}~~Since $h_M^{-p}: S^{n-1}\rightarrow\mathbb{R}$ is a bounded Borel function, let $f=gh_M^{-p}$ in (\ref{3.4}), we have
$$\int_{S^{n-1}}g(v)h_M^{-p}(v)d\widetilde{\mathcal{C}}_{q,j}(M,Q,v)
=\frac{1}{n}\int_{S^{n-1}}g(\alpha_M(u))h_M^{-p}(\alpha_M(u))\rho_M^q(u)\rho_Q^{n-q-j}(u)du,$$
which, in light of (\ref{3.3}), is the desired result (\ref{3.7}).      \hfill${\square}$

\noindent{\bf Lemma 3.3}~~{\it Suppose $p,q\in\mathbb{R}$ and $j\neq n$. If $M\in\mathcal{K}_o^n$ and $Q\in\mathcal{S}_o^n$, then for Borel set $\eta\subseteq S^{n-1}$ and each bounded, Borel function $g: S^{n-1}\rightarrow\mathbb{R}$
\begin{align}\label{3.8}
\nonumber&\int_{S^{n-1}}g(v)d\widetilde{\mathcal{C}}_{p,q,j}(M,Q,v)\\
&=\frac{1}{n}\int_{\partial'M}g(\nu_M(x))(x\cdot\nu_M(x))^{1-p}|x|^{-j}\|x\|_Q^{q+j-n}d\mathcal{H}^{n-1}(x),
\end{align}
\begin{align}\label{3.9}
\widetilde{\mathcal{C}}_{p,q,j}(M,Q,\eta)
=\frac{1}{n}\int_{x\in\chi}(x\cdot\nu_M(x))^{1-p}|x|^{-j}\|x\|_Q^{q+j-n}d\mathcal{H}^{n-1}(x),
\end{align}
where $\chi=\mathbb{X}_M(\eta)$.}

{\it \bf Proof.}~~Let $f=g(\alpha_M)h_M^{-p}(\alpha_M)\rho_Q^{n-q-j}(u)$ in (\ref{2.12}). From the homogeneity of $\rho_Q$, (\ref{2.5}), (\ref{2.6}) and finally the fact that $\|\cdot\|_Q=\rho_Q^{-1}$, combining with (\ref{3.7}), we get
\begin{align*}
&\int_{S^{n-1}}g(v)d\widetilde{\mathcal{C}}_{p,q,j}(M,Q,v)\\
&=\frac{1}{n}\int_{S^{n-1}}g(\alpha_M(u))h_M^{-p}(\alpha_M(u))\rho_M^q(u)\rho_Q^{n-q-j}(u)du\\
&=\frac{1}{n}\int_{\partial'M}g(\alpha_M(\overline{x}))h_M^{-p}(\alpha_M(\overline{x}))
|x|^{q-n}\rho_Q^{n-q-j}(\overline{x})(x\cdot\nu_M(x))d\mathcal{H}^{n-1}(x)\\
&=\frac{1}{n}\int_{\partial'M}g(\nu_M(x))(x\cdot\nu_M(x))^{1-p}|x|^{-j}\|x\|^{q+j-n}_Qd\mathcal{H}^{n-1}(x).
\end{align*}
This gives the (\ref{3.8}).

Let $g=\mathbb{I}_{\eta}$ in (\ref{3.8}). Recall that $\nu_M(x)\in\eta\Leftrightarrow x\in \mathbb{X}_M(\eta)$, for almost all $x$ with respect to spherical Lebesgue measure. Thus, we establish (\ref{3.9}).   \hfill${\square}$

\noindent{\bf Proposition 3.1}~~{\it For $p,q\in\mathbb{R}$ and $j\neq n$. If $M\in\mathcal{K}_o^n$ and $Q\in\mathcal{S}_o^n$, then
\begin{align}\label{3.10}
\widetilde{\mathcal{C}}_{p,q,0}(M,M,\cdot)
=\widetilde{\mathcal{C}}_{p,n,0}(M,B,\cdot)=\frac{1}{n}S_{p}(M,\cdot),
\end{align}
\begin{align}\label{3.11}
\widetilde{\mathcal{C}}_{0,q,j}(M,Q,\cdot)=\widetilde{\mathcal{C}}_{q,j}(M,Q,\cdot).
\end{align}}

{\it \bf Proof.}~~For all $x\in\partial M$, using the (\ref{3.1}), we have $\|x\|_M=\rho_M^{-1}(x)=1$ from the definition of the radial function. From (\ref{3.9}), we have
\begin{align*}
\widetilde{\mathcal{C}}_{p,n,0}(M,B,\eta)&=\frac{1}{n}\int_{x\in\chi}(x\cdot\nu_M(x))^{1-p}d\mathcal{H}^{n-1}(x)\\
&=\frac{1}{n}\int_{x\in\chi}(x\cdot\nu_M(x))^{1-p}\|x\|_M^{q+j-n}d\mathcal{H}^{n-1}(x)\\
&=\widetilde{\mathcal{C}}_{p,q,0}(M,M,\eta),\ \ \chi=\mathbb{X}_M(\eta).
\end{align*}
The (\ref{2.9}) states that the integral above is just $\frac{1}{n}S_{p}(M,\cdot)$, then (\ref{3.10}) and (\ref{3.11}) are obtained. From the (\ref{3.6}) and let $p=0$ in (\ref{3.6}), we get the (\ref{3.11}).  \hfill${\square}$

\section{\bf The $(p,q)$-mixed quermassintegrals}

$\quad$In this section, in order to prove the main results, we first give some variational formulas for $q$-th dual mixed quermassintegrals of a logarithmic family of Wulff shapes and convex hull. Then, we establish the integral representation and study some properties of $(p,q)$-mixed quermassintegrals.

Suppose $\Omega$ be a closed subset of $S^{n-1}$ that is not contained in any closed hemisphere. Let $f: \Omega\rightarrow\mathbb{R}$ be continuous and $\delta>0$. Define a positive continuous function, $h_t:\Omega\rightarrow(0,\infty)$, by letting
$$\log h_t(v)=\log h_0(v)+tf(v)+o(t,v),$$
for each $t\in(-\delta,\delta)$, where $o(t,\cdot): \Omega\rightarrow\mathbb{R}$ is continuous and $\lim_{t\rightarrow0}o(t,\cdot)/t=0$, uniformly on $\Omega$. Denote by
$$[h_t]=\{x\in\mathbb{R}: x\cdot v\leq h_t(v)\ for \ all \ v\in\Omega\},$$
the Wulff shape determined by $h_t$. We shall call $[h_t]$ a logarithmic family of Wulff shapes generated by $(h_0,f)$. If $h_0$ is the support function $h_M$ of a convex body $M$, we also write $[h_t]$ as $[M,f,t]$.

The Wulff shape of a function and the convex hull generated by its reciprocal are relate (see\cite{H}).

\noindent{\bf Lemma 4.1}~~{\it Suppose $\Omega\subset S^{n-1}$ is a closed set not contained in any closed hemisphere of $S^{n-1}$. Let $h: \Omega\rightarrow(0,\infty)$ be continuous. Then the Wulff shape $[h]$ determined by $h$ and the convex hull $\langle1\backslash h\rangle$ generated by the function $1\backslash h$ are polar reciprocals of each other; i.e.,
$$[h]^\ast=\langle1\backslash h\rangle.$$}

The following lemma shows that the support function of convex hull are differentiable.

\noindent{\bf Lemma 4.2} \cite{H}~~{\it Suppose $\Omega\subset S^{n-1}$ is a closed set not contained in any closed hemisphere of $S^{n-1}$. Let $\rho_0: \Omega\rightarrow(0,\infty)$ and $f: \Omega\rightarrow\mathbb{R}$ are continuous. If $\langle \rho_t\rangle$ is a logarithmic family of convex hulls generated by $(\rho_0,f)$, then, for $q\in\mathbb{R}$
$$\lim_{t\rightarrow0}\frac{h^{-q}_{\langle\rho_t\rangle}(v)-h^{-q}_{\langle\rho_0\rangle}(v)}{t}
=-qh^{-q}_{\langle\rho_0\rangle}(v)f(\alpha^\ast_{\langle\rho_0\rangle}(v)).$$
for all $v\in S^{n-1}\backslash\eta_{\langle\rho_0\rangle}$.}

The following theorem shows a variational formula for $q$-th dual mixed quermassintegrals of $q$-th dual mixed curvature measure and polar convex hull.

\noindent{\bf Theorem 4.1}~~{\it Suppose $\Omega\subset S^{n-1}$ is a closed set not contained in any closed hemisphere of $S^{n-1}$. Let $\rho_0: \Omega\rightarrow(0,\infty)$ and $f: \Omega\rightarrow\mathbb{R}$ are continuous. If $\langle\rho_t\rangle$ is a logarithmic family of convex hulls of $(\rho_0,f)$, and $q\neq0$, $j\neq n$, then
$$\lim_{t\rightarrow0}\frac{\widetilde{W}_{q,j}(\langle\rho_t\rangle^\ast,Q)
-\widetilde{W}_{q,j}(\langle\rho_0\rangle^\ast,Q)}{t}
=-q\int_{\Omega}f(u)d\widetilde{\mathcal{C}}_{q,j}(\langle\rho_0\rangle^\ast,Q,u)$$
which holds for $Q\in\mathcal{S}_o^n$.}

{\it \bf Proof of Theorem 4.1.}~~First, we write $\eta_{\langle\rho_0\rangle}$ by $\eta_0$. Recall that $\eta_0$ is the set of spherical Lebesgue measure zero that consists of the complement, in $S^{n-1}$, of the regular normal vectors of the convex body $\langle\rho_0\rangle=conv\{\rho_0(u)u: u\in\Omega\}$. The continuous function
$$\alpha_{\langle\rho_0\rangle}^\ast: S^{n-1}\backslash\eta_0\rightarrow S^{n-1}$$
is defined by $\alpha_{\langle\rho_0\rangle}^\ast(v)\in\Re_{\langle\rho_0\rangle}^\ast(v)=\{\alpha_{\langle\rho_0\rangle}^\ast(v)\}$ for all $v\in S^{n-1}\backslash\eta_0$

For $v\in S^{n-1}\backslash\eta_0$. Now, we prove that $\Re_{\langle\rho_0\rangle}(v)\subset\Omega$. Let
$$h_{\langle\rho_0\rangle}(v)=\max_{u\in\Omega}\rho_0(u)u\cdot v=\rho_0(u_0)u_0\cdot v$$
for some $u_0\in\Omega$. This means that
$$\rho_0(u_0)u_0\in H_{\langle\rho_0\rangle}(v),$$
thus $\rho_0(u_0)u_0\in \partial\langle\rho_0\rangle$ because in addition to $\rho_0(u_0)u_0$ obviously belonging to $\langle\rho_0\rangle$, it also belongs to $H_{\langle\rho_0\rangle}(v)$. Since $v$ is a regular normal vector of $\langle\rho_0\rangle$, therefore $\alpha_{\langle\rho_0\rangle}^\ast(v)=u_0\in\Omega$. Thus
\begin{align}\label{4.1}
\Re^\ast_{\langle\rho_0\rangle}(S^{n-1}\backslash\eta_0)\subset\Omega.
\end{align}
From (\ref{4.1}) and Lemma 2.1, we get the fact that
\begin{align}\label{4.2}
\Re_{\langle\rho_0\rangle^\ast}(S^{n-1}\backslash\eta_0)\subset\Omega.
\end{align}

Since $\Omega$ is closed, using the Tietze extension theorem, we can extend the continuous function $f: \Omega\rightarrow\mathbb{R}$ to a continuous function $\hat{f}: S^{n-1}\rightarrow\mathbb{R}$. From (\ref{4.2}), we get
\begin{align}\label{4.3}
f(\alpha_{\langle\rho_0\rangle^\ast}(v))=(\hat{f}\mathbb{I}_\Omega)(\alpha_{\langle\rho_0\rangle^\ast}(v)).
\end{align}

From (\ref{1.1}), (\ref{2.1}), Lemma 4.2, Lemma 2.1, (\ref{4.3}), (\ref{3.4}), and dominated convergence theorem. Since $\eta_0$ has measures zero, we get for $Q\in\mathcal{S}_o^n$
\begin{align*}
&\lim_{t\rightarrow0}\frac{\widetilde{W}_{q,j}(\langle\rho_t\rangle^\ast,Q)
-\widetilde{W}_{q,j}(\langle\rho_0\rangle^\ast,Q)}{t}\\
&=\lim_{t\rightarrow0}\frac{1}{n}\int_{S^{n-1}}\frac{\rho_{\langle\rho_t\rangle^\ast}^q(v)
-\rho_{\langle\rho_0\rangle^\ast}^q(v)}{t}\rho_{Q}^{n-q-j}(v)dv\\
&=\frac{1}{n}\int_{S^{n-1}\backslash\eta_0}\lim_{t\rightarrow0}\frac{h_{\langle\rho_t\rangle}^{-q}(v)
-h_{\langle\rho_0\rangle}^{-q}(v)}{t}\rho_{Q}^{n-q-j}(v)dv\\
&={-\frac{q}{n}}\int_{S^{n-1}\backslash\eta_0}f(\alpha^\ast_{\langle\rho_0\rangle}(v))
h^{-q}_{\langle\rho_0\rangle}(v)\rho_{Q}^{n-q-j}(v)dv\\
&={-\frac{q}{n}}\int_{S^{n-1}\backslash\eta_0}f(\alpha_{\langle\rho_0\rangle^\ast}(v))
\rho^{q}_{\langle\rho_0\rangle^\ast}(v)\rho_{Q}^{n-q-j}(v)dv\\
&={-\frac{q}{n}}\int_{S^{n-1}\backslash\eta_0}(\hat{f}\mathbb{I}_\Omega)(\alpha_{\langle\rho_0\rangle^\ast}(v))
\rho^{q}_{\langle\rho_0\rangle^\ast}(v)\rho_{Q}^{n-q-j}(v)dv\\
&=-q\int_{S^{n-1}}(\hat{f}\mathbb{I}_\Omega)(u)d\widetilde{\mathcal{C}}_{q,j}(\langle\rho_0\rangle^\ast,Q,u)\\
&=-q\int_{\Omega}f(u)d\widetilde{\mathcal{C}}_{q,j}(\langle\rho_0\rangle^\ast,Q,u).
\end{align*}
This give the proof of Theorem 4.1.    \hfill${\square}$

The following theorem shows a variational formula for $q$-th dual mixed quermassintegrals of $q$-th dual mixed curvature measure and Wulff shape.

\noindent{\bf Theorem 4.2}~~{\it Suppose $\Omega\subset S^{n-1}$ is a closed set not contained in any closed hemisphere of $S^{n-1}$. If $h_0: \Omega\rightarrow(0,\infty)$ and $f: \Omega\rightarrow\mathbb{R}$ are continuous, and $[h_t]$ is a logarithmic family of Wulff shapes generated by $(h_0,f)$, then
$$\lim_{t\rightarrow0}\frac{\widetilde{W}_{q,j}([h_t],Q)-\widetilde{W}_{q,j}([h_0],Q)}{t}
=q\int_{\Omega}f(v)d\widetilde{\mathcal{C}}_{q,j}([h_0],Q,v)$$
which holds for $Q\in\mathcal{S}_o^n$, $j\neq q\neq 0$.}

{\it \bf Proof of Theorem 4.2.}~~The logarithmic family of Wulff shapes $[h_t]$ is defined as the Wulff shape of $h_t$, where $h_t$ if given by
$$\log h_t(v)=\log h_0(v)+tf(v)+o(t,v).$$
Let $\rho_t=h_t^{-1}$, we get
$$\log \rho_t(v)=\log \rho_0(v)-tf(v)-o(t,v).$$
Let $\langle\rho_t\rangle$ be the logarithmic family of convex hulls associated with $(\rho_0,-f)$. From Lemma 4.1, we see that
$$[h_t]=\langle\rho_t\rangle^\ast,$$
and the desired conclusions now follows from Theorem 4.1.   \hfill${\square}$

Then we state the special case of Theorem 4.2 for logarithmic families of Wulff shapes generated by convex bodies. Here, we shall write $[h_t]$ as $[h,f,t]$, and if $h$ happens to be the support function of a convex body $M$ perhaps as $[M,f,t]$.

\noindent{\bf Theorem 4.3}~~{\it Suppose $M\in\mathcal{K}_o^n$ and $f: S^{n-1}\rightarrow\mathbb{R}$ are continuous. If $Q\in\mathcal{S}_o^n$, $j\neq n$ and $q\neq 0$, then
$$\lim_{t\rightarrow0}\frac{\widetilde{W}_{q,j}([M,f,t],Q)-\widetilde{W}_{q,j}(M,Q)}{t}
=q\int_{S^{n-1}}f(v)d\widetilde{\mathcal{C}}_{q,j}(M,Q,v).$$}

From (\ref{3.7}), let $Q=M$ and $g=h_N^p$, we see that
$$\int_{S^{n-1}}g(v)d\widetilde{\mathcal{C}}_{p,q,j}(M,M,v)
=\frac{1}{n}\int_{S^{n-1}}h_N^p(\alpha_M(u))h_M^{-p}(\alpha_M(u))\rho_M^{n-j}(u)du.$$
According to the (\ref{3.10}) and (\ref{2.10}), the $L_p$ mixed quermassintegrals, $W_{p,j}(M,N)$, has an integral formula: For $M,N\in\mathcal{K}_o^n$, then
\begin{align}\label{4.4}
W_{p,j}(M,N)=\frac{1}{n}\int_{S^{n-1}}\bigg(\frac{h_N}{h_M}\bigg)^p(\alpha_M(u))\rho_M^{n-j}(u)du.
\end{align}

Using the above results, we obtain the following important result, namely, the variational formulas of $q$-th dual mixed quermassintegrals with respect to $(p,q)$-dual mixed curvature measures.

\noindent{\bf Theorem 4.4}~~{\it Suppose $p,q\neq0$, $j\neq n$. If $M,N\in\mathcal{K}_o^n$ and $Q\in\mathcal{S}_o^n$, then
\begin{align}\label{4.5}
\lim_{t\rightarrow0}\frac{\widetilde{W}_{q,j}(M+_pt\cdot N,Q)-\widetilde{W}_{q,j}(M,Q)}{t}
=\frac{q}{p}\int_{n-1}h_N^p(v)d\widetilde{C}_{p,q,j}(M,Q,v).
\end{align}}

{\it \bf Proof of Theorem 4.4.}~~Choose small enough $\varepsilon$, define $h_t$ by
\begin{align}\label{4.6}
h_t^p=h_M^p+t h_N^p,\ \ p\neq0.
\end{align}
From (\ref{2.2}) and (\ref{4.6}), the Wulff shape $[h_t]=M+_pt\cdot N$. Let $f=\frac{1}{p}\frac{h_N^p}{h_M^p}$. By Theorem 4.3, we have
\begin{align}\label{4.7}
\lim_{t\rightarrow0}\frac{\widetilde{W}_{q,j}(M+_pt\cdot N,Q)-\widetilde{W}_{q,j}(M,Q)}{t} =\frac{q}{p}\int_{S^{n-1}}h_N^p(v)h_M^{-p}(v)d\widetilde{\mathcal{\mathcal{C}}}_{q,j}(M,Q,v).
\end{align}
From (\ref{3.3}), using (\ref{3.7}) with $g=h_N^p$, we can write (\ref{4.7}) as
$$\lim_{t\rightarrow0}\frac{\widetilde{W}_{q,j}(M+_pt\cdot N,Q)-\widetilde{W}_{q,j}(M,Q)}{t}=\frac{q}{p}\int_{S^{n-1}}h_N^p(v)d\widetilde{\mathcal{C}}_{p,q,j}(M,Q,v)$$
Hence, above result complete the proof of Theorem 4.4.   \hfill${\square}$

\noindent{\bf Corollary 4.1}~~{\it Let $j=0$ in (4.5), the result of $\widetilde{V}_{p,q}(M,N,Q)=\widetilde{W}_{p,q,0}(M,N,Q)$ was obtained in \cite{L4}.}

In order to obtain a unification which includes the $L_p$ mixed quermassintegrals and dual mixed quermassintegrals, this leads us to define the $(p,q)$-dual mixed quermassintegrals.

\noindent{\bf Definition 4.1}~~{\it Suppose $p,q\in\mathbb{R}$, $j\neq n$. If $M,N\in\mathcal{K}_o^n$ and $Q\in\mathcal{S}_o^n$, define the $(p,q)$-dual mixed quermassintegrals, $\widetilde{W}_{p,q,j}(M,N,Q)$, by
\begin{align}\label{4.7.1}
\frac{q}{p}\widetilde{W}_{p,q,j}(M,N,Q)=\lim_{t\rightarrow0}\frac{\widetilde{W}_{q,j}(M+_pt\cdot N,Q)-\widetilde{W}_{q,j}(M,Q)}{t}
\end{align}}

By using (\ref{3.7}) with $g=h_N^p$, the Definition 4.1 can be written as a following formula,
\begin{align}\label{4.7.2}
\widetilde{W}_{p,q,j}(M,N,Q)=\frac{1}{n}\bigg(\frac{h_N}{h_M}\bigg)^p(\alpha_M(u))\bigg(\frac{\rho_M}{\rho_Q}\bigg)^q(u)\rho_Q^{n-j}(u)du.
\end{align}

Now, we study some properties of $(p,q)$-mixed quermassintegrals.

\noindent{\bf Proposition 4.1}~~{\it Suppose $p,q\in\mathbb{R}$, $j\neq n$. If $M,N\in\mathcal{K}_o^n$ and $Q\in\mathcal{S}_o^n$, then
\begin{align}\label{4.8}
\widetilde{W}_{p,q,j}(M,M,M)=\widetilde{W}_{j}(M),
\end{align}
\begin{align}\label{4.9}
\widetilde{W}_{p,q,j}(M,M,Q)=\widetilde{W}_{q,j}(M,Q),
\end{align}
\begin{align}\label{4.10}
\widetilde{W}_{p,q,j}(M,N,M)=W_{p,j}(M,N),
\end{align}
\begin{align}\label{4.11}
\widetilde{W}_{0,q,j}(M,N,Q)=\widetilde{W}_{q,j}(M,Q),
\end{align}
\begin{align}\label{4.12}
\widetilde{W}_{p,n,j}(M,N,Q)=W_{p,j}(M,N),
\end{align}
\begin{align}\label{4.13}
\widetilde{W}_{p,q,0}(M,N,Q)=\widetilde{V}_{p,q}(M,N,Q).
\end{align}}

{\it \bf Proof.}~~From (\ref{4.7.2}) and the polar coordinate formula for volume, the (\ref{4.8}) is obtained. By (\ref{4.7.2}) and the definition of dual mixed quermassintegrals (\ref{1.1}), it yields (\ref{4.9}). Together with (\ref{4.7.2}) and the definition of mixed quermassintegrals (\ref{4.4}), identit (\ref{4.10}) is obtained. Identity (\ref{4.11}) follow from (\ref{4.7.2}) and the definition of dual mixed quermassintegrals (\ref{1.1}). Note that (\ref{4.7.2}) and the definition of mixed quermassintegrals (\ref{4.4}) it easily conclude that (\ref{4.12}) is ture. Finally, when $j=0$ in (\ref{4.7.2}), then $(p,q)$-mixed volume is obtained.    \hfill${\square}$

Next, we give the property of general linear transformation under $GL(n)$, the group of general linear transformation.

\noindent{\bf Proposition 4.2}~~{\it Support $p,q\in\mathbb{R}$, $j\neq n$. If $M,N\in\mathcal{K}_o^n$ and $Q\in\mathcal{S}_o^n$, then
\begin{align}\label{4.14}
\widetilde{W}_{p,q,j}(\phi M,\phi N,\phi Q)=|det\phi|\widetilde{W}_{p,q,j}(M,N,Q),
\end{align}
for each $\phi\in GL(n)$.}

{\it \bf Proof.}~~For $\phi\in GL(n)$. From the definition (\ref{4.7.1}) and (\ref{2.3}), combined with the fact that $\widetilde{W}_{q,j}(\phi M,\phi Q)=|det\phi|\widetilde{W}_{q,j}(M,Q)$, we get
\begin{align*}
&\frac{q}{p}\widetilde{W}_{p,q,j}(\phi M,\phi N,\phi Q)\\
&=\lim_{t\rightarrow0}\frac{\widetilde{W}_{q,j}(\phi M+_pt\cdot (\phi N),\phi Q)-\widetilde{W}_{q,j}(\phi M,\phi Q)}{t}\\
&=\lim_{t\rightarrow0}\frac{\widetilde{W}_{q,j}(\phi(M+_pt\cdot N),\phi Q)-\widetilde{W}_{q,j}(\phi M,\phi Q)}{t}\\
&=|det\phi|\lim_{t\rightarrow0}\frac{\widetilde{W}_{q,j}(M+_pt\cdot N,Q)-\widetilde{W}_{q,j}(M,Q)}{t}\\
&=\frac{q}{p}|det\phi|\widetilde{W}_{p,q,j}(\phi M,\phi N,\phi Q),
\end{align*}
i.e.
$$\widetilde{W}_{p,q,j}(\phi M,\phi N,\phi Q)=|det\phi|\widetilde{W}_{p,q,j}(M,N,Q).$$
This give the (\ref{4.14}).    \hfill${\square}$

\noindent{\bf Corollary 4.2}~~If each $\phi\in SL(n)$ in (\ref{4.14}), then
$$\widetilde{W}_{p,q,j}(\phi M,\phi N,\phi Q)=\widetilde{W}_{p,q,j}(M,N,Q).$$

\noindent{\bf Proposition 4.3}~~{\it The $(p,q)$-mixed quermassintegrals $\widetilde{W}_{p,q,j}(M,N,Q)$ is a valuation respect to $N\in\mathcal{K}_o^n$, and is a valuation with respect to $Q\in\mathcal{S}_o^n$.}

{\it \bf Proof.}~~From Definition 4.1
\begin{align}\label{4.15}
\widetilde{W}_{p,q,j}(M,N,Q)
=\frac{1}{n}\int_{S^{n-1}}h_N^p(\alpha_M(u))h_M^{-p}(\alpha_M(u))\rho_M^q(u)\rho_Q^{n-q-j}(u)du.
\end{align}
The $(p,q)$-mixed quermassintegrals is a valuation on $\mathcal{S}_o^n$ with respect to the $Q$ can be seen easily by writing as
and observing that for $Q_1,Q_2\in\mathcal{S}_o^n$, we get
$$\rho_{Q_1\cup Q_2}^{n-q-j}+\rho_{Q_1\cap Q_2}^{n-q-j}=\rho_{Q_1}^{n-q-j}+\rho_{Q_2}^{n-q-j}.$$

That the $(p,q)$-mixed quermassintegrals is a valuation on $\mathcal{K}_o^n$ with respect to the $N$ can be see easily by  looking at (\ref{4.15}) and using the fact that if $N_1,N_2\in\mathcal{K}_o^n$, are such that $N_1\cup N_2\in\mathcal{K}_o^n$, we get
$$h_{N_1\cup N_2}^p+h_{N_1\cap N_2}^p=h_{N_1}^p+h_{N_2}^p.$$  \hfill${\square}$

\section{\bf Geometric inequalities}

~~~~In this section, we obtain some important geometric inequalities of $(p,q)$-mixed quermassintegrals.  Firstly, we give the $L_p$ Minkowski type inequality for $(p,q)$-mixed quermassintegrals.

\noindent{\bf Theorem 5.1}~~{\it Suppose $p,q\in\mathbb{R}$ are such that $1\leq\frac{q}{n-j}\leq p$, $j\neq n$. If $M,N\in\mathcal{K}_o^n$ and $Q\in\mathcal{S}_o^n$, then
\begin{align}\label{5.1}
\widetilde{W}_{p,q,j}(M,N,Q)\geq W_j(M)^{\frac{q-p}{n-j}}W_j(N)^{\frac{p}{n-j}}\widetilde{W}_j(Q)^{\frac{n-q-j}{n-j}},
\end{align}
with equality if and only if $M,N,Q$ are dilates when $1<\frac{q}{n-j}<p$, while only $M$ and $N$ are dilates when $q=n$ and $p>1$, and $M$ and $N$ are homotheic when $q=n$ and $p=1$.}

{\it \bf Proof of Theorem 5.1.}~~From (\ref{4.7.2}), (\ref{1.2}) and Definition 4.1, if $1\leq\frac{q}{n-j}\leq p$, by the H\"{o}lder inequality, we have
\begin{align*}
&\widetilde{W}_{p,q,j}(M,N,Q)
=\frac{1}{n}\int_{S^{n-1}}\bigg(\frac{h_N}{h_M}\bigg)^p(\alpha_M(u))\rho_M^q(u)\rho_Q^{n-q-j}(u)du\\
&=\frac{1}{n}\int_{S^{n-1}}\bigg[\bigg(\frac{h_N}{h_M}\bigg)^{\frac{(n-j)p}{q}}
(\alpha_M(u))\rho_M^{n-j}(u)\bigg]^{\frac{q}{n-j}}
\bigg[\rho_Q^{n-j}(u)\bigg]^{\frac{n-q-j}{n-j}}du\\
&\geq\bigg[\frac{1}{n}\int_{S^{n-1}}\bigg(\frac{h_N}{h_M}\bigg)^{\frac{(n-j)p}{q}}(\alpha_M(u))
\rho_M^{n-j}(u)du\bigg]^{\frac{q}{n-j}}
\bigg[\frac{1}{n}\int_{S^{n-1}}\rho_Q^{n-j}(u)du\bigg]^{\frac{n-q-j}{n-j}}\\
&=W_{\frac{(n-j)p}{q},j}(M,N)^{\frac{q}{n-j}}\widetilde{W}_j(Q)^{\frac{n-q-j}{n-j}}\\
&\geq W_j(M)^{\frac{q-p}{n-j}}W_j(N)^{\frac{p}{n-j}}\widetilde{W}_j(Q)^{\frac{n-q-j}{n-j}}.
\end{align*}
The equality condition of (\ref{5.1}) follow from the equality conditions of H\"{o}lder inequality and (\ref{2.11}).    \hfill${\square}$

\noindent{\bf Corollary 5.1}~~{\it Let $j=0$ in Theorem 5.1, we can get the Minkowski inequality of $(p,q)$-mixed volumes (see e.g., \cite{L4}).}

Then, we show the monotonic inequality of $(p,q)$-mixed quermassintegrals.

\noindent{\bf Theorem 5.2}~~{\it Suppose $p,q\in\mathbb{R}$ are such that $0<p<q$ or $p<0< q$. If $M,N\in\mathcal{K}_o^n$, $Q\in\mathcal{S}_o^n$ and $j\neq n$, then
\begin{align}\label{5.2}
\bigg(\frac{\widetilde{W}_{p,q,j}(M,N,Q)}{\widetilde{W}_{q,j}(M,Q)}\bigg)^{\frac{1}{p}}
\geq\bigg(\frac{\widetilde{W}_{p-q,q,j}(M,N,Q)}{\widetilde{W}_{q,j}(M,Q)}\bigg)^{\frac{1}{p-q}},
\end{align}
with equality if and only if $M,N,Q$ are dilates.}

{\it \bf Proof of Theorem 5.2.}~~From (\ref{4.7.2}), H\"{o}lder inequality, we obtain that for $\frac{p}{p-q}>1$,
\begin{align*}
&\widetilde{W}_{p,q,j}(M,N,Q)
=\frac{1}{n}\int_{S^{n-1}}\bigg(\frac{h_N}{h_M}\bigg)^p(\alpha_M(u))\bigg(\frac{\rho_M}{\rho_Q}\bigg)^q(u)
\rho_Q^{n-j}(u)du\\
&=\frac{1}{n}\int_{S^{n-1}}\bigg[\bigg(\frac{h_N}{h_M}\bigg)^{p-q}(\alpha_M(u))\bigg(\frac{\rho_M}{\rho_Q}\bigg)^q(u)
\rho_Q^{n-j}(u)\bigg]^{\frac{p}{p-q}}\bigg[\bigg(\frac{\rho_M}{\rho_Q}\bigg)^q(u)\rho_Q^{n-j}(u)\bigg]^{\frac{-q}{p-q}}du\\
&\geq\bigg[\frac{1}{n}\int_{S^{n-1}}\bigg(\frac{h_N}{h_M}\bigg)^{p-q}(\alpha_M(u))\bigg(\frac{\rho_M}{\rho_Q}\bigg)^q(u)
\rho_Q^{n-j}(u)du\bigg]^{\frac{p}{p-q}}\\
&\ \ \ \ \cdot\bigg[\frac{1}{n}\int_{S^{n-1}}\bigg(\frac{\rho_M}{\rho_Q}\bigg)^q(u)\rho_Q^{n-j}(u)du\bigg]^{\frac{-q}{p-q}}\\
&=\widetilde{W}_{p-q,q,j}(M,N,Q)^{\frac{p}{p-q}}\widetilde{W}_{q,j}(M,Q)^{-\frac{q}{p-q}}.
\end{align*}
This implies that for $p>0$
$$\bigg(\frac{\widetilde{W}_{p,q,j}(M,N,Q)}{\widetilde{W}_{q,j}(M,Q)}\bigg)^{\frac{1}{p}}
\geq\bigg(\frac{\widetilde{W}_{p-q,q,j}(M,N,Q)}{\widetilde{W}_{q,j}(M,Q)}\bigg)^{\frac{1}{p-q}}.$$
According to equality condition of H\"{o}lder inequality, we see that equality holds in (\ref{5.2}) if and only if $M,N,Q$ are dilates.    \hfill${\square}$

\noindent{\bf Corollary 5.2}~~{\it Let $\widetilde{W}_j(Q)$ replace $\widetilde{W}_{q,j}(M,Q)$ in Theorem 1.3, we can get $$\bigg(\frac{\widetilde{W}_{n-p,q,j}(M,N,Q)}{\widetilde{W}_{j}(M)}\bigg)^{\frac{1}{n-p}}
\geq\bigg(\frac{\widetilde{W}_{n-q,q,j}(M,N,Q)}{\widetilde{W}_{j}(M)}\bigg)^{\frac{1}{n-q}},$$
with equality if and only if $M,Q$ are dilates.}

\noindent{\bf Corollary 5.3}~~{\it Let $j=0$ in Corollary 5.2, we can get the monotonic inequality of $(p,q)$-mixed volumes (see e.g., \cite{FH}).}

Finally, we obtain a type of cyclic inequalities of $(p,q)$-mixed quermassintegrals as follows.

\noindent{\bf Theorem 5.3}~~{\it Suppose $p,q,r,s\in\mathbb{R}$ satisfy $ p<q<r\leq n$. If $M,N\in\mathcal{K}_o^n$, $Q\in\mathcal{S}_o^n$ and $j\neq n$, then
\begin{align}\label{5.3}
\widetilde{W}_{q,s,j}(M,N,Q)^{r-p}\leq\widetilde{W}_{p,s,j}(M,N,Q)^{r-q}\widetilde{W}_{r,s,j}(M,N,Q)^{q-p},
\end{align}
with equality if and only if $M,N$ and $Q$ are dilates.}

\noindent{\bf Theorem 5.4}~~{\it Suppose $p,q,r,s\in\mathbb{R}$ satisfy $p<q<r\leq n$. If $M,N\in\mathcal{K}_o^n$, $Q\in\mathcal{S}_o^n$ and $j\neq n$, then
\begin{align}\label{5.4}
\widetilde{W}_{s,q,j}(M,N,Q)^{r-p}\leq\widetilde{W}_{s,p,j}(M,N,Q)^{r-q}\widetilde{W}_{s.r,j}(M,N,Q)^{q-p},
\end{align}
with equality if and only if $M,N$ and $Q$ are dilates.}

{\it \bf Proof of Theorem 5.3.}~~Suppose $p,q,r,s\in\mathbb{R}$ satisfy $p<q<r\leq n$. From (\ref{4.7.2}) and H\"{o}lder inequality, we get that for $u\in S^{n-1}$
\begin{align*}
&\widetilde{W}_{p,s,j}(M,N,Q)^{\frac{r-q}{r-p}}\widetilde{W}_{r,s,j}(M,N,Q)^{\frac{q-p}{r-p}}\\
&=\bigg[\frac{1}{n}\int_{S^{n-1}}\bigg(\frac{h_N}{h_M}\bigg)^p(\alpha_M(u))\rho_M^{s}(u)
\rho_Q^{n-s-j}(u)du\bigg]^{\frac{r-q}{r-p}}\\
&\ \ \ \ \cdot\bigg[\frac{1}{n}\int_{S^{n-1}}\bigg(\frac{h_N}{h_M}\bigg)^r(\alpha_M(u))\rho_M^{s}(u)
\rho_Q^{n-s-j}(u)du\bigg]^{\frac{q-p}{r-p}}\\
&=\bigg[\frac{1}{n}\int_{S^{n-1}}\bigg(\bigg(\frac{h_N}{h_M}\bigg)^{\frac{p(r-q)}{r-p}}(\alpha_M(u))
\rho_M^{\frac{s(r-q)}{r-p}}(u)\rho_Q^{\frac{(n-s-j)(r-q)}{r-p}}(u)\bigg)^{\frac{r-p}{r-q}}du\bigg]^{\frac{r-q}{r-p}}\\
&\ \ \ \ \cdot\bigg[\frac{1}{n}\int_{S^{n-1}}\bigg(\bigg(\frac{h_N}{h_M}\bigg)^{\frac{r(q-p)}{r-p}}(\alpha_M(u))
\rho_M^{\frac{s(q-p)}{r-p}}(u)\rho_Q^{\frac{(n-s-j)(q-p)}{r-p}}(u)\bigg)^{\frac{r-p}{q-p}}du\bigg]^{\frac{q-p}{r-p}}\\
&\geq\frac{1}{n}\int_{S^{n-1}}\bigg(\frac{h_N}{h_M}\bigg)^q(\alpha_M(u))\rho_M^{s}(u)\rho_Q^{n-s-j}(u)du\\
&=\widetilde{W}_{q,s,j}(M,N,Q).
\end{align*}
i.e.
$$\widetilde{W}_{q,s,j}(M,N,Q)^{r-p}\leq\widetilde{W}_{p,s,j}(M,N,Q)^{r-q}\widetilde{W}_{r,s,j}(M,N,Q)^{q-p}.$$
This yields (\ref{5.3}). According to the equality condition of H\"{o}lder's integral inequality, we see that equality holds in (\ref{5.3}) if and only if $M,N$, and $Q$ are dilates.    \hfill${\square}$

{\it \bf Proof of Theorem 5.4.}~~Similar to the proof of Theorem 5.3, we can easily get Theorem 5.4.  \hfill${\square}$

\noindent{\bf Corollary 5.4}~~{\it Let $j=0$ in Theorem 5.3 and 5.4, we can get the cyclic inequalities of $(p,q)$-mixed volumes (see e.g., \cite{C}).}

\noindent{\bf Corollary 5.5}\cite{L3}~~{\it Let $Q=M$ in Theorem 5.3, and together with (4.10), we can get the monotonic inequality of $L_p$-mixed quermassintegrals
$$W_{q,j}(M,N)^{r-p}\geq W_{p,j}(M,N)^{r-q}W_{r,j}(M,N)^{q-p},$$
with equality if and only if $M$ and $N$ are dilates.

If $N=M$ in Theorem 5.4, by (\ref{4.9}), we can get another monotonic inequality of dual mixed quermassintegrals
$$\widetilde{W}_{q,j}(M,N)^{r-p}\geq \widetilde{W}_{p,j}(M,N)^{r-q}\widetilde{W}_{r,j}(M,N)^{q-p},$$
with equality if and only if $M$ and $N$ are dilates.}

\section{\bf The $(p,q)$-dual mixed Minkowski problem}

~~~~The Minkowski problem and dual Minkowski problem are the hot topic of research in Brunn-Minkowski theory and dual Brunn-Minkowski theory. In this section, we further study the Minkowski problem for the $(p,q)$-dual mixed curvature measures.

Here, in order to prove the existence of Minkowski problem for $(p,q)$-dual mixed curvature measures, we can transform the existence problem into a maximization problem. Therefore, we need to define an important function below, which is the key to contact Minkowski problem and maximization problem.

Suppose $\mu$ is a non-zero finite Borel measure on $S^{n-1}$, $p,q\neq0$ and $j\neq n$, we define $\Phi_{p,q,j}: \mathcal{K}_o^n\rightarrow\mathbb{R}$ by
\begin{align}\label{6.1}
\Phi_{p,q,j}(M,Q)
=-\frac{1}{p}\log\int_{S^{n-1}}h_M^p(v)d\mu(v)+\frac{1}{q}\log\int_{S^{n-1}}h_M^q(u)\rho_Q^{n-q-j}(u)du,
\end{align}
for any $M\in\mathcal{K}_o^n$, and $Q\in\mathcal{S}_o^n$.

\noindent{\bf Theorem 6.1}~~{\it (Existence) For $p,q\in\mathbb{R}$, $j\neq n$ and $Q\in\mathcal{S}_o^n$. If $\mu$ is a Borel measure on $S^{n-1}$, then there exists a convex body $M\in\mathcal{K}_o^n$ such that $\widetilde{\mathcal{C}}_{p,q,j}(M,Q,\cdot)=\mu$.}

Next, we consider the maximization problem.
$$\sup\{\Phi_{p,q,j}(M,Q):M\in\mathcal{K}_o^n\}.$$

In the following lemma, we can prove that a solution to the Minkowski problem for measure $\mu$ is also a solution to a maximization problem for the function $\Phi_{p,q,j}$.

\noindent{\bf Lemma 6.1}~~{\it Suppose $p,q\neq0$, $j\neq n$ and $\mu$ is a non-zero finite Borel measure on $S^{n-1}$. For fixed $Q\in\mathcal{S}^n_o$, if $M\in\mathcal{K}_o^n$ and satisfies the following conditions
\begin{align}\label{6.2}
\int_{S^{n-1}}h_M^p(v)d\mu(v)=\widetilde{W}_{q,j}(M,Q),
\end{align}
and
\begin{align}\label{6.3}
\Phi_{p,q,j}(M,Q)=\sup\{\Phi_{p,q,j}(N,Q): N\in\mathcal{K}_o^n\},
\end{align}
then
$$\mu=\widetilde{\mathcal{C}}_{p,q,j}(M,Q,\cdot),$$}

{\it \bf Proof.}~~For any $f,g\in C^+(S^{n-1})$, corresponding to (\ref{6.1}), we first define the functional $\Psi_{p,q,j}: C^+(S^{n-1})\rightarrow\mathbb{R}$ by
$$\Psi_{p,q,j}(f,g)=-\frac{1}{p}\log\int_{S^{n-1}}f^p(v)d\mu(v)
+\frac{1}{q}\log\int_{S^{n-1}}\rho_{[f]}^q(u)\rho_{\langle g\rangle}^{n-q-j}(u)du,$$
where $[f]$ denotes the Wulff shape generated by $f$ and $\langle g\rangle$ denotes the convex hull generated by $g$.

To prove the final result, we claim that
\begin{align}\label{6.4}
\Psi_{p,q,j}(f,g)\leq \Psi_{p,q,j}(h_M,\rho_Q),
\end{align}
for each $f,g\in C^+(S^{n-1})$ and $M\in\mathcal{K}_o^n$, $Q\in\mathcal{S}^n_o$.

We give the following facts to explain our claim. From the notions of Wullf shape and convex hull, we have $h_{[f]}\leq f$, $\rho_{\langle g\rangle}\geq g$, $[h_M]=M$ and $\langle \rho_Q\rangle=Q$. Thus
$$\Psi_{p,q,j}(f,g)\leq \Psi_{p,q,j}(h_{[f]},\rho_{\langle g\rangle})=\Phi_{p,q,j}([f],\langle g\rangle)
\leq \Phi_{p,q,j}(N,Q)=\Psi_{p,q,j}(h_M,\rho_Q),$$
where the second inequality sign follows from (\ref{6.3}).

For any $z\in C(S^{n-1})$ and $t\in(-\delta,\delta)$ where $\delta>0$ is sufficiently small, let
\begin{align}\label{6.5}
h_t(v)=h_M(v)e^{tz(v)}.
\end{align}
From Theorem 4.2, we have
\begin{align}\label{6.6}
\frac{d}{dt}\widetilde{W}_{q,j}([h_t],Q)\bigg|_{t=0}=q\int_{S^{n-1}}z(v)d\widetilde{\mathcal{C}}_{q,j}(N,Q).
\end{align}
According to the (\ref{6.4}), the definition of $\Psi_{p,q,j}$ and (\ref{6.6}), we have
\begin{align*}
0&=\frac{d}{dt}\Psi_{p,q,j}(h_t,g)\\
&=\frac{d}{dt}\bigg(-\frac{1}{p}\log\int_{S^{n-1}}h_t^p(v)d\mu(v)
+\frac{1}{q}\log\widetilde{W}_{q,j}([h_t],\langle\rho_Q\rangle)\bigg)\bigg|_{t=0}\\
&=-\bigg(\int_{S^{n-1}}h^p_M(v)d\mu(v)\bigg)^{-1}\int_{S^{n-1}}h_M^p(v)z(v)d\mu(v)\\
&\ \ \ \ +\frac{1}{\widetilde{W}_{q,j}(M,Q)}\int_{S^{n-1}}z(v)d\widetilde{\mathcal{C}}_{q,j}(M,Q,v).
\end{align*}
Using (\ref{6.2}), then
$$h_N^p(v)d\mu(v)=d\widetilde{\mathcal{C}}_{q,j}(M,Q,v).$$
This together with (\ref{3.3}), we get
$\mu=\widetilde{\mathcal{C}}_{p,q,j}(M,Q,\cdot)$.  \hfill${\square}$

Therefore, we have shown that the Minkowski problem can be transformed into a maximization problem via Lemma 6.1. The next key task is to prove the existence of solutions to the maximum problem. Before we do that, let us give the following lemmas. The first lemma was established in \cite{HZ}.

\noindent{\bf Lemma 6.2}~~{\it For $p,\varepsilon_0>0$ and $\mu$ is a non-zero even finite Borel measure on $S^{n-1}$. Suppose $e_{1k},...,e_{nk}$ is a sequence of orthonormal basis in $\mathbb{R}^n$ and $\{a_k\}$ is a sequence of positive real numbers. Assume $e_{1k},...,e_{nk}$ converges to an orthonormal basis $e_{1},...,e_{n}$ in $\mathbb{R}^n$. Define
$$G_k=\{x\in\mathbb{R}^n: |x\cdot e_{1k}|^2+\cdot\cdot\cdot+x\cdot e_{n-1,k}|^2\leq a_{k}^2\ \
and \ \ x\cdot e_{nk}|^2\leq\varepsilon_0\}.$$
If $\mu$ is not concentrated in any great subsphere, then there exists $c,L>0$ (independent of $k$) such that
$$\int_{S^{n-1}}h_{G_k}^p(v)d\mu(v)\geq c,$$
for any $k>L$.}

The following lemma estimates the $q$-th dual mixed quermassintegrals of the generalized ellipsoid $T_\alpha$.

\noindent{\bf Lemma 6.3}~~{\it Suppose $0<\alpha<1$ and $e_1,...,e_n$ be an orthonormal basis in $\mathbb{R}^n$. If $0<q<n-j$, $j\neq n$ and fixed $Q\in\mathcal{S}_o^n$, then
$$\lim_{\alpha\rightarrow0^+}\int_{S^{n-1}}\rho_{T_\alpha}^q(u)\rho_Q^{n-q-j}(u)du=0,$$
where $T_\alpha$ is defined by
$$T_{\alpha}=\{x\in\mathbb{R}^n: |x\cdot e_1|\leq\alpha\ \ and \ \ |x\cdot e_2|^2+\cdot\cdot\cdot|x\cdot e_n|^2\leq1\},$$
for $0<\alpha<1$ and $e_1,...,e_n$ be the orthonormal basis in $\mathbb{R}^n$.}

{\it \bf Proof.}~~For $0<\alpha<1$ and fixed $Q\in\mathcal{S}_o^n$, then $\rho_{T_\alpha}$ and $\rho_Q$ are bounded. Therefore, using the dominated convergence theorem, we get
$$\lim_{\alpha\rightarrow0^+}\int_{S^{n-1}}\rho_{T_\alpha}^q(u)\rho_Q^{n-q-j}(u)du
=\int_{S^{n-1}}\lim_{\alpha\rightarrow0^+}\rho_Q^{n-q-j}(u)\rho_{T_\alpha}^q(u)du.$$
It is easy to see that
$$\lim_{\alpha\rightarrow0^+}\rho_Q^{n-q-j}(u)\rho_{T_\alpha}^q(u)=\bigg\{_{0,\ \ \ \ \ \ \ \ \ \ \ \ \ \  otherwise.}^{\rho_Q^{n-q-j}(u),\ \ \ \ u\in span\{e_1,...,e_n\}}$$
Therefore,
$$\lim_{\alpha\rightarrow0^+}\int_{S^{n-1}}\rho_{T_\alpha}^q(u)\rho_Q^{n-q-j}(u)du=0.$$   \hfill${\square}$

Now, we establish the existence of maximization problem.

\noindent{\bf Lemma 6.4}~~{\it Suppose $p,q>0, p\neq q, j\neq n$ and $\mu$ is a non-zero even finite Borel measure on $S^{n-1}$. If $\mu$ is not concentrated in any great subsphere, then there exists $M\in\mathcal{K}_e^n$ such that
\begin{align}\label{6.7}
\Phi_{p,q,j}(M,Q)=\sup\{\Phi_{p,q,j}(N,Q): N\in\mathcal{K}_e^n\},
\end{align}
for fixed $Q\in\mathcal{S}_o^n$.}

Let $\{N_k\}\subset\mathcal{K}_e^n$ be a maximizing sequence; i.e.,
$$\lim_{k\rightarrow\infty}\Phi_{p,q,j}(N_k,Q)=\sup\{\Phi_{p,q,j}(N,Q): N\in\mathcal{K}_e^n\},$$
for fixed $Q\in\mathcal{S}_o^n$.

Combining with Lemma 6.2 and Lemma 6.3, for $Q=B$, Lemma 6.4 was proved in \cite{HZ}. When $Q$ is a star body in $\mathcal{S}_o^n$, the proof method of Lemma 6.4 is similar, and thus will be omitted.

{\it \bf Proof of Theorem 6.1.}~~The proof follows directly from the Lemma 6.1 and Lemma 6.4.   \hfill${\square}$

In order to prove the uniqueness for the $(p,q)$-dual mixed Minkowski problem, we need the following facts.

{\it $(p,q)$-dual mixed curvature measures of polytopes.} For $j\neq n$ and $M\in\mathcal{K}_o^n$ be a polytope with outer unit normals $v_1,...,v_m$. Let $\triangle_j$ be the cone that consists of the set of all rays emanating from the origin and passing through the facet of $M$ whose outer normal is $v_j$. Recalling that we abbreviate $\Re_M^\ast(\{v_j\})$ by $\Re_M^\ast(v_j)$, from the definition of reverse radial Gauss image, we get
\begin{align}\label{6.8}
\Re_M^\ast(v_j)=S^{n-1}\cap\triangle_j,\ \ and \ \alpha_M(u)=v_j,\ for\ almost\ all\ u\in S^{n-1}\cap\triangle_j.
\end{align}
If $\eta\subset S^{n-1}$ is a Borel set such that $\{v_1,...,v_m\}\cap\eta=\varnothing$, then $\Re_M^\ast(\eta)$ has spherical Lebesgue measure zero. Thus, the $(p,q)$-dual mixed curvature measure $\widetilde{\mathcal{C}}_{p,q,j}(M,Q,\cdot)$ is discrete and is concentrated on $\{v_1,...,v_m\}$. From the (\ref{3.6}), and (\ref{6.8}), we get
\begin{align}\label{6.9}
\widetilde{\mathcal{C}}_{p,q,j}(M,Q,\cdot)
=\frac{1}{n}\sum_{j=1}^m \delta_{v_j}h_M^{-p}(v_j)\int_{S^{n-1}\cap\triangle_j}\rho_M^q(u)\rho_Q^{n-q-j}(u)du,
\end{align}
where $\delta_{v_j}$ denotes the delta measure concentrated at $v_j$.

\noindent{\bf Theorem 6.2}~~{\it (Uniqueness) For $p,q\in\mathbb{R}$, $j\neq n$, $Q\in\mathcal{S}_o^n$ and $M_1,M_2\in\mathcal{K}_o^n$ be polytopes. If
$$\widetilde{\mathcal{C}}_{p,q,j}(M_1,Q,\cdot)=\widetilde{\mathcal{C}}_{p,q,j}(M_2,Q,\cdot),$$
then $M_1=M_2$ when $p<q$, and $M_1$ is a dilate of $M_2$ when $q=p$.}

{\it \bf Proof of Theorem 6.2.}~~From the first fact, we see that the $(p,q)$-dual mixed curvature measures of polytopes are discrete. If $\widetilde{\mathcal{C}}_{p,q,j}(M_1,Q,\cdot)=\widetilde{\mathcal{C}}_{p,q,j}(M_2,Q,\cdot)$. By (\ref{6.9}), then $M_1$ and $M_2$ have the same outer unit normals, we have
\begin{align}\label{6.10}
\nonumber&\frac{1}{n}h_{M_1}^{-p}(v_j)\int_{S^{n-1}\cap\triangle_j}\rho_{M_1}^q(u)\rho_Q^{n-q-j}(u)du\\
&=\frac{1}{n}h_{M_2}^{-p}(v_j)\int_{S^{n-1}\cap\triangle_j^{\prime}}\rho_{M_2}^q(u)\rho_Q^{n-q-j}(u)du,
\end{align}
where $\triangle_j$ and $\triangle_j^{\prime}$ are the cones formed by the origin and the facets of $M_1$ and $M_2$ with normal $v_j$, respectively.

If $M_1\neq M_2$, clearly $M_1\subsetneq M_2$ is not possible. Let $\lambda$ be the maximal number such that $\lambda M_1\subseteq M_2$, then $\lambda<1$. Since $\lambda M_1$ and $M_2$ have the same outer unit normals, there is a facet of $\lambda M_1$ that is contained in a facet of $M_2$.

Let $j=0$ in (\ref{6.10}), then the proof of uniqueness was obtain in \cite{L4}. If $j\neq n$, the proof of uniqueness for the $(p,q)$-dual mixed Minkowski problem is very similar, therefore, we can get the result.       \hfill${\square}$

\end{document}